\documentclass[12pt]{amsart}

\newtheorem{theorem}{Theorem}[section]
\newtheorem{proposition}[theorem]{Proposition}
\newtheorem{lemma}[theorem]{Lemma}
\newtheorem{corollary}[theorem]{Corollary}
\theoremstyle{definition}   
\newtheorem{definition}{Definition}

\theoremstyle{remark}

\numberwithin{equation}{section}

\usepackage{times}
\usepackage{enumerate}
\usepackage{tfrupee}
\usepackage{tikz}
\usetikzlibrary{chains,fit}
\usetikzlibrary{shapes,snakes}
\usetikzlibrary{graphs}
\usepackage{graphicx,adjustbox}
\usepackage{hyperref}
\usepackage{amsmath,amssymb}
\usetikzlibrary{cd}

\newcommand{\lcm}{{\rm \mbox{lcm}}}

\newcommand{\supp}{\mbox{\rm supp}}
\newcommand{\mdeg}{{\rm \mbox{mdeg}}}

\usepackage{amscd}
\usepackage{graphicx}
\usepackage[all,cmtip]{xy}

\title[Transversal Intersection of Monomial Ideals]
{Transversal Intersection of Monomial Ideals}
\author{
Joydip Saha
\and
Indranath Sengupta
\and
Gaurab Tripathi
}
\date{}

\address{\small \rm  Discipline of Mathematics, IIT Gandhinagar, Palaj, Gandhinagar, 
Gujarat 382355, INDIA.} 
\email{saha.joydip56@gmail.com}
\thanks{The first author is a post-doctoral research fellow under the 
research project EMR/2015/000776 sponsored by the SERB, Government of India.}

\address{\small \rm  Discipline of Mathematics, IIT Gandhinagar, Palaj, Gandhinagar, 
Gujarat 382355, INDIA.}
\email{indranathsg@iitgn.ac.in}
\thanks{The second author is the corresponding author, who is supported by the 
research project EMR/2015/000776 sponsored by the SERB, Government of India.}

\address{\small \rm Department of Mathematics, Jadavpur University, Kolkata,
WB 700 032, India.} 
\email{gelatinx@gmail.com}

\date{}



\begin{document}
\begin{abstract}
In this paper we prove conditions for transversal intersection 
of monomial ideals and derive a simplicial characterization of this 
phenomenon.
\end{abstract}

\maketitle

\section{Introduction}
The transversal intersection phenomenon for monomial ideals 
in the polynomial ring happens to be extremely interesting 
and it turns out that it is equivalent to having disjoint 
supports for their minimal generating sets; see Theorem 
\ref{disjoint}. The Taylor complex which resolves (possibly 
non-minimally) a monomial ideal has been understood completely 
for ideals of the form $I+J$, when $I$ and $J$ are monomial 
ideals intersecting transversally; see Theorem \ref{tcomplex}. 
As an application of this theorem we prove in Corollary \ref{tenminimal} 
that for two ideals $I$ and $J$ in $R$ intersecting transversally, 
the ideal $I+J$ is resolved minimally by the complex $\mathbb{M}(I)_{\centerdot}\otimes \mathbb{M}(J)_{\centerdot}$, if $\mathbb{M}(I)_{\centerdot}$ and $\mathbb{M}(J)_{\centerdot}$ denote minimal 
free resolutions of $I$ and $J$ respectively. Minimal free resolutions 
for ideals of the form $I+J$ have an interesting structure when $I$ 
and $J$ intersect transversally and are supported simplicially; see \ref{simpsupport}. 

\section{Monomial Ideals}
Let $R=K[x_{1},x_{2},\cdots ,x_{n}]$, where $x_{i}$'s are indeterminates over 
the field $K$. Let ${\rm Mon}(R)$ denote the set of all monomials in $R$. 
Every nonzero polynomial $f\in R$ is a unique $K$-linear combination 
of monomials given by $f=\sum_{v\in {\rm Mon}(R)} a_{v}v$. Let 
$m(f):= \{v\in {\rm Mon}(R) \mid a_{v}\neq 0\}$. An ideal 
$I$ in $R$ is said to be a \textit{monomial ideal} if it is 
generated by monomials of $R$. We list down some standard facts on 
monomial ideals. 
\medskip

\begin{proposition}\label{mingenmono}
\begin{enumerate}
\item $I$ is a monomial ideal if and only if $f\in  I \Rightarrow m(f)\subset I$.
\medskip

\item Let $\lbrace u_{1},\cdots ,u_{m}\rbrace$ be a monomial generating set of an 
ideal $I$, where $u_{i}$'s are monomials. A monomial $v\in I$ if and only if 
$v=u_{i}w$, for some $1\leq i\leq m$.
\medskip

\item Each monomial ideal $I$ has a unique minimal monomial set of generators $G(I)$.
\end{enumerate}
\end{proposition}

\proof See  \cite{miller}.\qed
\medskip

\begin{definition}Let $\emptyset\neq T\subset {\rm Mon}(R)$. We define 
$$
{\rm supp}(T) = \{i \mid x_{i}\ {\rm divides}\ m \ \text{for\, some}\ m\in T\}.
$$
If $T = \{m\}$, we simply write ${\rm supp}(m)$ instead of ${\rm supp}(\{m\})$. 
If $S$ and $T$ are two nonempty subsets of ${\rm Mon}(R)$, then, 
${\rm supp}(S)\cap {\rm supp}(T)=\emptyset$ if and only if 
${\rm supp}(f)\cap {\rm supp}(g)=\emptyset$ for every $f\in S$ and $g\in T$.
\end{definition}
\medskip

\begin{definition}
We say that the ideals $I$ and $J$ of $R$ intersect transversally 
if $I\cap J= IJ$.
\end{definition}
\medskip

\begin{theorem}\label{disjoint}
Let $I$ and $J$ be two monomial ideals of $R$. Then, $I\cap J =IJ$ if and only if 
$\supp(G(I)) \cap \supp(G(J))=\emptyset$.
\end{theorem}
\proof Let $ I\cap J =IJ$. Consider the set 
$$S=\{{\rm lcm}(f,g)\mid f\in G(I)\,, g\in G(J) \quad {\rm and}\quad \supp(f)\cap \supp(g)\neq 
\emptyset\}.$$
If $S= \emptyset$, then for every $f \in G(I)$ and $g \in G(J)$ we have 
$\supp(f)\cap \supp(g)=\emptyset$, which proves that 
$\supp(G(I)) \cap \supp(G(J))=\emptyset$.
\medskip

Suppose that $S\neq \emptyset$. We define a partial order $\leq$ on $S$ in the 
following way: Given $s,t\in S$, we define $s\leq t$ if and only 
$s\mid t$. Let $\mathbf{m}$ be a minimal element (by Zorn's lemma) of $S$; then there 
exists $f,g$ such that $f\in G(I)$ and $g\in G(J)$ such that $\lcm(f,g)=\mathbf{m}$. 
Since $\lcm(f,g)\in I\cap J$ and $I\cap J=IJ$, we have $\lcm(f,g)=\mathbf{m}\in IJ$. 
By Proposition \ref{mingenmono}, there exist $h_{1}\in G(I)$ and $h_{2}\in G(J)$, 
such that $h_{1}h_{2}\mid \mathbf{m}$, since the generating set of $IJ $ is the set 
$G(I)G(J)=\lbrace uv\mid u\in G(I), v\in G(J)\rbrace$. 
\medskip

We prove that $\supp(h_{1})\cap \supp(h_{2}) = \emptyset$. If this is not the case, then 
$\lcm(h_{1}, h_{2})\neq  h_{1}h_{2}$, in fact $\lcm(h_{1}, h_{2})< h_{1}h_{2}$ and 
$\lcm(h_{1}, h_{2})\mid \mathbf{m}$; contradicting minimality of $\mathbf{m}$ in $S$. 
\medskip

We now prove that $\supp(f)\cap \supp(h_{2})\neq \emptyset$. We know that 
$h_{1}h_{2}\mid \lcm(f,g)=\mathbf{m}$. Therefore, if 
$\supp(f)\cap \supp(h_{2})= \emptyset$, then $h_{2}\mid g$; which contradicts 
minimality of the generating set $G(J)$. Similarly we can prove that 
$\supp(g)\cap \supp(h_{1})\neq \emptyset$. Now, $h_{1}h_{2}\mid \mathbf{m}=\lcm(f,g)$ 
implies that $h_{2}\mid \mathbf{m}=\lcm(f,g)$. Moreover, $f\mid \lcm(f,g)$. Therefore, 
$\lcm(f,h_{2})\mid \lcm(f,g)$. Similarly, $\lcm(g,h_{1})\mid \lcm(f,g)$. Now if 
$\lcm(f,h_{2})< \mathbf{m}$ or $\lcm(g,h_{1})<\lcm(f,g)=\mathbf{m}$, then we 
have a contradiction, since $\lcm(f,g)$ and $\lcm(g,h_{1})$ both are in the set $S$ 
and $\mathbf{m}$ is a minimal element in $S$. Therefore, we have must have 
$\lcm(f,h_{2})=\lcm(g,h_{1})=\lcm(f,g)=\mathbf{m}=h_{1}h_{2}w$ (say) and 
$\dfrac{fh_{2}}{\rm{gcd}(f,h_{2})}=\dfrac{gh_{1}}{\gcd(g,h_{1})}=h_{1}h_{2}w$, 
therefore $\dfrac{f}{\gcd(f,h_{2})}=h_{1}w$. Hence, $h_{1}\mid f$ and that 
contradicts minimality of the generating set $G(I)$. Hence, 
$S=\emptyset$ and we are done.
\medskip

Conversely, let us assume that $\supp(G(I))\cap\supp(G(J)) = \emptyset$. 
Without loss of generality, we can assume that $\supp(G(I)) = \{1,2,\ldots , k\}$ and 
$\supp(G(J)) = \{k+1, k+2,..., n\}$. Let $f\in I\cap J$, such that 
$f=\sum_{v\in Mon(R)} a_{v}v$. We have $v\in I\cap J$ for all $v\in m(f)$. 
It is therefore enough to show that if $m$ is a monomial and 
$m\in I\cap J$, then $m\in IJ$. Let $m\in I\cap J$; there exist 
$m_{I}\in G(I)$ and $m_{J}\in G(J)$ such that $m_{I}\mid m$ and $m_{J}\mid m$. 
Let $m_{I}m_{1}=m$; then $m_{J}\mid m=m_{I}m_{1}$. We know that 
$\supp(m_{I})\cap \supp(m_{J}) = \emptyset$, therefore $m_{J}\mid m_{1}$ 
and it follows that $m_{1}=m_{J}m_{2}$ and $m=m_{I}m_{J}m_{2}$. Therefore, 
$m\in IJ$ since $m_{I}m_{J}\in IJ$. \qed
\medskip

\begin{theorem}\label{dissup}
Let $I$ and $J$ be two ideals of $R$ such that w.r.t. some monomial order $\preceq$, 
$\supp(\rm{Lt}(I)) \cap \supp(\rm{Lt}(J))=\emptyset$. Then $I$ and $J$ intersect 
transversally.
\end{theorem}

\proof Let $\mathcal{G}_{I}$ and $\mathcal{G}_{J}$ be Gr\"{o}bner bases with respect 
to the monomial order $\preceq$ of the ideals $I$ and $J$ respectively. Again we assume 
that $f\in I\cap J$, then there exist polynomials $p\in \mathcal{G}_{I}$ and 
$q\in \mathcal{G}_{J}$ such that $\rm{Lt}(p)\mid \rm{Lt}(f)$ and 
$\rm{Lt}(q)\mid \rm{Lt}(f)$. Since $\supp(\rm{Lt}(I)) \cap \supp(\rm{Lt}(J))=\emptyset$, 
we have  $\rm{Lt}(p)\cdot\rm{Lt}(q)\mid \rm{Lt}(f)$. After division we write $f=pq+r$, 
then $r\in I\cap J$ and $\rm{Lt}(r)\preceq\rm{Lt}(f)$. we can apply same process on $r$, 
and after finite stage we get $f\in IJ$. \qed 
\bigskip

\section{The Taylor complex} 
We define the {\em multidegree} of a 
monomial $x^{\mathbf{a}} = x_{1}^{a_{1}}x_{2}^{a_{2}}\cdots x_{n}^{a_{n}}$ to be the 
$n$-tuple $(a_{1}, \ldots , a_{n})$, denoted by 
$\text{mdeg}(\mathbf{x}^{\mathbf{a}})$. Hence $R=\oplus _{m\in {\rm Mon}(R)} Km$ 
has a direct sum decomposition as $K$-vector spaces $Km$, where $Km = \{cm\mid c\in K\}$.
\medskip

\begin{definition}[\textbf{The Taylor complex}] 
Let $M$ be a monomial ideal of $R$ minimally generated by the 
monomials $m_{1},m_{2},\cdots , m_{p}$. The Taylor complex 
for $M$ is given by $\mathbb{T}(M)_{\centerdot}$ as follows: $\mathbb{T}(M)_{i}$ is the 
free $R$-module generated by $e_{j_{1}}\wedge e_{j_{2}}\wedge \cdots\wedge e_{j_{i}}$, for all $1\leq j_{1} < \cdots < j_{i}\leq p$, where 
$\{e_{1}, \ldots , e_{p}\}$ is the standard basis for the free $R$-module $R^{p}$. 
The differential 
$\delta(M)_{\cdot}$ is given by the following:
\medskip

\noindent $\delta(M)_{i}(e_{j_{1}}\wedge \cdots\wedge e_{j_{i}})$
$$=\sum_{1\leq k\leq i}(-1)^{k-1}
\dfrac{\lcm (m_{j_{1}}, \cdots ,m_{j_{i}})}{\lcm(m_{j_{1}}, \cdots,\widehat{m_{j_{k}}}, \cdots, m_{j_{i}})}(e_{j_{1}}\wedge \cdots\wedge \widehat{e_{j_{k}}}\wedge\cdots\wedge e_{j_{i}}).$$
We define ${\rm mdeg}(e_{j_{1}}\wedge \cdots\wedge e_{j_{i}})
= \mdeg(\lcm (m_{j_{1}}\cdots m_{j_{i}}))$. The Taylor complex 
$(\mathbb{T}(M)_{\centerdot}, \delta(M)_{\centerdot})$ defined above is in deed a chain complex 
of free $R$ modules and gives a free resolution for the $R$-module $R/M$; 
see section 26 in \cite{peeva}.

\end{definition}
\medskip

\begin{theorem}The Taylor complex $(\mathbb{T}(M)_{\centerdot}, \delta(M)_{\centerdot})$ 
defined above is a free resolution (though not minimal) of the monomial ideal M; 
called the Taylor resolution of $M$.
\end{theorem}  

\proof See Theorem 26.7 in \cite{peeva}.\qed
\medskip

\begin{theorem}\label{tcomplex}
Let $I$ and $J$ be two monomial ideals of $R$ intersecting transversally. Let 
$\mathbb{T}(I)_{\centerdot}$ and $\mathbb{T}(J)_{\centerdot}$ be the Taylor resolutions 
of $I$ and $J$ respectively. Then, the Taylor's resolution of $I+J$ is isomorphic 
to the complex $\mathbb{T}(I)_{\centerdot}\otimes \mathbb{T}(J)_{\centerdot}$ and hence 
$\mathbb{T}(I)_{\centerdot}\otimes \mathbb{T}(J)_{\centerdot}$ is acyclic. 
\end{theorem}

\proof We know that $\mathbb{T}(I)_{i}= R^{\binom{p}{i}}$ and $\mathbb{T}(J)_{j} = R^{\binom{q}{j}}$, where $p$ and $q$ denote the minimal number of generators 
of $I$ and $J$ respectively. Then, 
\begin{eqnarray*}
(\mathbb{T}(I)_{\centerdot}\otimes \mathbb{T}(J)_{\centerdot})_{r} & = &\oplus_{i+j=r}(\mathbb{T}(I)_{i}\otimes \mathbb{T}(J)_{j})\\
{} &  = & \oplus_{i+j=r} R^{\binom{p}{i}}\otimes R^{\binom{q}{j}}\\
{} &  = & R^{\sum_{i+j=r}\binom{p}{i}\binom{q}{j}}\\
{} & = & R^{\binom{p+q}{r}} = \mathbb{T}(I+J)_{r}.
\end{eqnarray*}

A basis of $(\mathbb{T}(I)_{\centerdot}\otimes \mathbb{T}(J)_{\centerdot})_{r}$ is given by $(e_{k_{1}}\wedge \cdots\wedge e_{k_{i}})\otimes (e'_{k'_{1}}\wedge\cdots\wedge e'_{k'_{j}})$, 
such that $i+j=r$ and $1\leq k_{1} < k_{2} < \cdots < k_{i}\leq p$ and 
$1\leq k'_{1} < k'_{2} < \cdots < k'_{j}\leq q$. Now for each $r$, the free 
module $(\mathbb{T}(I)_{\centerdot}\otimes \mathbb{T}(J)_{\centerdot})_{r}$ can be graded with the help of mdeg as follows. We define 
\begin{eqnarray*}
\mdeg\left((e_{k_{1}}\wedge \cdots\wedge e_{k_{i}})\otimes (e'_{k'_{1}}\wedge \cdots\wedge e'_{k'_{j}})\right)\\
= \mdeg\left(\lcm( m_{k_{1}},\ldots ,m_{k_{i}})\cdot \lcm(m'_{k'_{1}},\ldots , m'_{k'_{j}})\right).
\end{eqnarray*} 
This defines a multi-graded structure for the complex $\mathbb{T}(I)_{\centerdot}\otimes \mathbb{T}(J)_{\centerdot}$.
\medskip

Let $G(I)=\{m_{1},m_{2}\cdots,m_{p}\}$ and $G(J)=\{m'_{1},m'_{2}\cdots,m'_{q}\}$ 
denote the minimal set of generators for $I$ and $J$ respectively. Then $I+J$ is 
minimally generated by $G(I)\cup G(J)$, by Theorem \ref{disjoint}, since $G(I)$ 
and $G(J)$ are of disjoint support. We will now show that the tensor product complex 
$\mathbb{T}(I)_{\centerdot}\otimes \mathbb{T}(J)_{\centerdot}$is isomorphic to the Taylor 
resolution $\mathbb{T}(I+J)_{\centerdot}$. We define $\psi_{r}:(\mathbb{T}(I)_{\centerdot}\otimes \mathbb{T}(J)_{\centerdot})_{r}\longrightarrow \left(\mathbb{T}(I+J)_{\centerdot}\right)_{r}$ as  
$$(e_{k_{1}}\wedge\cdots\wedge e_{k_{i}})\otimes (e'_{k'_{1}}\wedge\cdots\wedge e'_{k'_{j}})\mapsto (e_{k_{1}}\wedge\cdots\wedge e_{k_{i}}\wedge e'_{k'_{1}}\wedge\cdots\wedge e'_{k'_{j}}).$$
Moreover, 
\begin{eqnarray*}
{} & \hspace*{-1in} \mdeg(e_{k_{1}}\wedge \cdots\wedge e_{k_{i}}\wedge e_{k'_{1}}\wedge\cdots\wedge e_{k'_{j}})\\
= & \hspace*{-1in} \mdeg\left(\lcm( m_{k_{1}},\ldots , m_{k_{s}}, m^{'}_{k'_{1}},\ldots, m^{'}_{k'_{s'}})\right)\\
= & \hspace*{-6pt} \left(\mdeg(\lcm(m_{k_{1}},\ldots , m_{k_{s}})), \mdeg(\lcm(m^{'}_{k'_{1}},\ldots , m^{'}_{k'_{s'}}))\right),
\end{eqnarray*}
since $G(I)$ and $G(J)$ are with disjoint supports. Hence, the map is a 
graded isomorphism between the free modules $(\mathbb{T}(I)_{\centerdot}\otimes \mathbb{T}(J)_{\centerdot})_{r}$ and 
$\left(\mathbb{T}(I+J)_{\centerdot}\right)_{r}$. 
We therefore have the following diagram of complexes:
\medskip

$\xymatrix{
\cdots \ar[r]  & (\mathbb{T}(I)_{\centerdot}\otimes \mathbb{T}(J)_{\centerdot})_{r}\ar[r]^{d_{r}}\ar[d]^{\psi_{r}}&(\mathbb{T}(I)_{\centerdot}\otimes \mathbb{T}(J)_{\centerdot})_{r-1}\ar[r]\ar[d]^{\psi_{r-1}}&\cdots\\
\cdots\ar[r] & \left(\mathbb{T}(I+J)_{\centerdot}\right)_{r}\ar[r]^{\delta(I+J)_{r}} & \left(\mathbb{T}(I+J)_{\centerdot}\right)_{r-1}\ar[r]& \cdots}$
\medskip

\noindent We now show that $\delta_{r} \circ \psi_{r}  =  \psi_{r-1} \circ d_{r}$. Let us take an 
arbitrary basis element of $(\mathbb{T}(I)_{\centerdot}\otimes \mathbb{T}(J)_{\centerdot})_{r}$, which is of the form 
$$(e_{k_{1}}\wedge e_{k_{2}}\wedge\cdots\wedge e_{k_{i}})\otimes 
(e_{k'_{1}}\wedge e_{k'_{2}}\wedge\cdots\wedge e_{k'_{j}}), \quad {\rm with} \quad i+j=r.$$ 
Then,
\begin{eqnarray*}
{} & (\psi_{r-1} \circ d_{r})\left((e_{k_{1}}\wedge\cdots\wedge e_{k_{i}})\otimes (e'_{k'_{1}}\wedge\cdots\wedge e'_{k'_{j}})\right)\\
= & \psi_{r-1}(\delta(I)_{i}(e_{k_{1}}\wedge\cdots\wedge e_{k_{i}})
\otimes (e'_{k'_{1}}\wedge\cdots\wedge e'_{k'_{j}})\\
{} & + \quad (-1)^{i}(e_{k_{1}}\wedge\cdots\wedge e_{k_{i}})\otimes \delta(J)_{j}(e'_{k'_{1}}\wedge\cdots\wedge e'_{k'_{j}})).
\end{eqnarray*}
Let us write 
\smallskip

$Q_{t} = \dfrac{\lcm(m_{k_{1}},\ldots , m_{k_{i}})}{\lcm(m_{k_{1}},\ldots , \widehat{m_{k_{t}}},\ldots , m_{k_{i}})}$, \, $1\leq t\leq i$; 
\smallskip

$V_{t} = \dfrac{\lcm(m_{k_{1}},\ldots , m_{k_{i}}, m'_{k'_{1}},\ldots, m'_{k'_{j}})}{\lcm(m_{k_{1}}, \ldots , \widehat{m_{k_{t}}}, \ldots , m_{k_{i}}, m'_{k'_{1}}, \ldots, m'_{k'_{j}})}$, 
 \, $1\leq t\leq i$; 
\smallskip

\noindent We know that $G(I)$ and $G(J)$ have disjoint supports, therefore, 
$$V_{t} = \dfrac{\lcm(m_{k_{1}},\ldots , m_{k_{i}})\cdot \lcm( m'_{k'_{1}},\ldots, m'_{k'_{j}})}{\lcm(m_{k_{1}}, \ldots , \widehat{m_{k_{t}}}, \ldots , m_{k_{i}})\cdot\lcm( m'_{k'_{1}}, \ldots, m'_{k'_{j}})}=Q_{t}\quad \forall \, 1\leq t\leq i.$$
Then
\begin{eqnarray*}
{} & \left(\delta(I)_{i}(e_{k_{1}}\wedge\cdots\wedge e_{k_{i}})\right)\otimes (e'_{k'_{1}}\wedge\cdots\wedge e'_{k'_{j}})\\
= & \sum_{t=1}^{i}(-1)^{t}Q_{t}\cdot (e_{k_{1}}\wedge\cdots\wedge \widehat{e_{k_{t}}}\cdots\wedge e_{k_{i}})\otimes (e'_{k'_{1}}\wedge\cdots\wedge e'_{k'_{j}});
\end{eqnarray*}
and 
\begin{eqnarray*}
{} & \psi_{r-1}\left(\delta(I)_{i}(e_{k_{1}}\wedge\cdots\wedge e_{k_{i}})\otimes (e'_{k'_{1}}\wedge\cdots\wedge e'_{k'_{j}})\right)\\
= & \hspace*{-6pt} \sum_{t=1}^{i}(-1)^{t}Q_{t}\cdot (e_{k_{1}}\wedge\cdots\wedge \widehat{e_{k_{t}}}\cdots\wedge e_{k_{i}})\otimes (e'_{k'_{1}}\wedge\cdots\wedge e'_{k'_{j}})\\
= & \hspace*{-6pt}\sum_{t=1}^{i}(-1)^{t}V_{t}\cdot (e_{k_{1}}\wedge\cdots\wedge \widehat{e_{k_{t}}}\wedge e_{k_{i}}\wedge e'_{k'_{1}}\wedge\cdots\wedge e'_{k'_{j}}).
\end{eqnarray*}  
The last equality in the above expression follows from the fact that $G(I)$ and $G(J)$ have disjoint supports. Similarly, it can be proved that
\begin{eqnarray*}
{} & \psi_{r-1}\left((e_{k_{1}}\wedge\cdots\wedge e_{k_{i}})\otimes(\delta(J)_{j}(e'_{k'_{1}}\wedge\cdots\wedge e'_{k'_{j}}))\right)\\
= & \hspace*{-6pt}\sum_{l=1}^{j}(-1)^{l}W_{l}\cdot (e_{k_{1}} \wedge\cdots\wedge e_{k_{i}}\wedge e^{'}_{k'_{1}}\wedge\cdots \widehat{e^{'}_{k'_{l}}}\cdots\wedge 
e^{'}_{k'_{j}});
\end{eqnarray*}
where 
$$W_{l}=\dfrac{\lcm(m_{k_{1}}, \ldots, m_{k_{i}}, m'_{k'_{1}},\ldots, m'_{k'_{j}})}{\lcm(m_{k_{1}}, \ldots,  m_{k_{i}}, m'_{k'_{1}},\ldots , \widehat{m'_{k'_{l}}},\ldots, m'_{k'_{j}})}, \, 
1\leq l\leq j.$$
Hence,
\begin{eqnarray*}
(*) & \hspace*{-11pt}(\psi_{r-1} \circ d_{r})\left((e_{k_{1}}\wedge\cdots\wedge e_{k_{i}})\otimes (e'_{k'_{1}}\wedge\cdots\wedge e'_{k'_{j}})\right)\\
= & \hspace*{-11pt} \sum_{t=1}^{i}(-1)^{t}V_{t}\cdot (e_{k_{1}}\wedge\cdots\wedge \widehat{e_{k_{t}}}\wedge e_{k_{i}}\wedge e'_{k'_{1}}\wedge\cdots\wedge e'_{k'_{j}})\\
+ & \hspace*{-11pt} \sum_{l=1}^{j}(-1)^{i+l}W_{l}\cdot (e_{k_{1}} \wedge\cdots\wedge e_{k_{i}}\wedge e^{'}_{k'_{1}}\wedge\cdots \widehat{e^{'}_{k'_{l}}}\cdots\wedge 
e^{'}_{k'_{j}})
\end{eqnarray*}
We now introduce $i+j$ new symbols $w_{s}$, such that 
$w_{s}=m_{k_{s}}$ for $s\leq i$ and $w_{s}=m'_{k'_{s-i}}$ for $s>i$. We 
also introduce another set of $i+j$ symbols $E_{s}$, such $E_{s}=e_{k_{s}}$ 
for $s\leq i$ and $E_{s}=e'_{k'_{s-i}}$ for $s>i$. Hence, the expression after 
the last equality in (*) can be written in a compact form as 
$$\sum_{k=1}^{i+j}(-1)^{k}\dfrac{\lcm(w_{1}, \ldots, w_{i+j})}{\lcm(w_{1},\ldots,\widehat{w_{k}},\ldots, w_{i+j})}(E_{1}\wedge\cdots\wedge \widehat{E_{k}}\cdots\wedge E_{i+j}).$$ 
Now 
\begin{eqnarray*}
{} & \hspace*{-1.25in} (\delta_{r} \circ \psi_{r})((e_{k_{1}}\wedge\cdots\wedge e_{k_{i}})\otimes (e'_{k'_{1}}\wedge\cdots\wedge e'_{k'_{j}}))\\ 
= & \hspace*{-2in} \delta_{r}(e_{k_{1}}\wedge\cdots\wedge e_{k_{i}}\wedge e'_{k'_{1}}\wedge\cdots\wedge e'_{k'_{j}})\\
= & \sum_{k=1}^{i+j}(-1)^{k}\dfrac{\lcm(w_{1}, \ldots, w_{i+j})}{\lcm(w_{1},\ldots, \widehat{w_{k}}, \ldots , w_{i+j})}(E_{1}\wedge\cdots\wedge \widehat{E_{k}}\cdots\wedge E_{i+j}).
\end{eqnarray*}
Hence, the diagram is commutative.\qed
\medskip

\begin{corollary}\label{tenminimal}
Let $I$ and $J$ be ideals in $R$ such that $IJ = I\cap J$. Let 
$\mathbb{M}(I)_{\centerdot}$ and $\mathbb{M}(J)_{\centerdot}$ denote minimal free resolutions 
of $I$ and $J$ respectively. Then, $\mathbb{M}(I)_{\centerdot}\otimes \mathbb{M}(J)_{\centerdot}$ 
is a minimal free resolution of the ideal $I+J$. 
\end{corollary}

\proof The minimal free resolutions $\mathbb{M}(I)_{\centerdot}$ and $\mathbb{M}(J)_{\centerdot}$ are direct summands of the Taylor resolutions $\mathbb{T}(I)_{\centerdot}$ and $\mathbb{T}(J)_{\centerdot}$ 
respectively. It follows that $(\mathbb{M}(I)_{\centerdot}\otimes \mathbb{M}(J)_{\centerdot})$ is a 
direct summand of $\mathbb{T}(I)_{\centerdot}\otimes \mathbb{T}(J)_{\centerdot} = \mathbb{T}(I+J)_{\centerdot}$. Hence, $(\mathbb{M}(I)_{\centerdot}\otimes \mathbb{M}(J)_{\centerdot})$ is a free resolution of 
$I+J$ and it is minimal since $\mathbb{M}(I)_{\centerdot}$ and $\mathbb{M}(J)_{\centerdot}$ are both minimal.\qed
\bigskip

\section{A Simplicial characterization of transversal intersection of monomial ideals}
We first introduce the basic definitions of a simplicial complex; see \cite{miller}.
\medskip

\begin{definition}
A \textit{simplicial complex} $\Delta$ on the vertex set $\{1,\ldots,m\}$ is 
a collection of subsets called faces or simplices, satisfying the following condition: 
$\sigma\in \Delta$ and $\tau\subset\sigma$ implies that $\tau\in \Delta$. 
A simplex $\sigma\in \Delta$ of cardinality $|\sigma|=r+1$ has dimension 
$r$ and it is called an \textit{$r$ face} of $\Delta$.
\end{definition}
\medskip

\begin{definition}
A \textit{facet} is a  maximal simplex of a simplicial complex $\Delta$. 
Let $\Gamma(\Delta)$ denote the set of all facets of $\Delta$. Then, the 
vertex set of $\Delta$ and the set $\Gamma(\Delta)$ determine $\Delta$ completely.
\end{definition}
\medskip

\begin{definition}
A \textit{standard  simplicial complex} of dimension $m-1$ on the vertex set 
$\{1,\ldots,m\}$ is the simplicial complex whose $\Gamma(\Delta)=\{\{1,\ldots,m\}\}$.
\end{definition}
\medskip

\begin{definition}
Let $\Delta_{1}$ and $\Delta_{2}$ be simplicial complexes on disjoint vertex sets. 
Let the vertext set of $\Delta_{1}$ be $\{1,\ldots,m\}$ and that of 
$\Delta_{2}$ be $\{m+1,\ldots,m+p\}$. The \textit{join of simplicial complexes} 
$\Delta_{1}$ and $\Delta_{2}$ is the simplicial complex 
$\Delta_{1}*\Delta_{2}$, whose vertex set is $\{1,\ldots,m+p\}$ and 
$$\Gamma(\Delta_{1}*\Delta_{2})=\{\sigma_{1}\cup\sigma_{2}\mid \sigma_{1}\in \Gamma(\Delta_{1}), \, \sigma_{2}\in \Gamma(\Delta_{2}) \}.$$ 
\end{definition}
\medskip

We now define a \textit{frame} as a complex of $K$-vector spaces with a fixed basis, which encodes the minimal free resolution of a monomial ideal; see \cite{peeva}.
\medskip

\begin{definition}
An $r-$frame $\mathbb{U}_{\centerdot}$ is  a complex of finite $K$ vector spaces with 
differential $\partial$ and a fixed basis satisfying the following:
\begin{enumerate}[(i)]
\item $\mathbb{U}_{i}=0$ for $i<0$;
\item $\mathbb{U}_{0}=K$;
\item $\mathbb{U}_{1}=K^{r}$, with basis $\{w_{1}, \ldots , w_{r}\}$;
\item $\partial(w_{j})=1$, for all $j = 1, \ldots , r$.
\end{enumerate}
\end{definition}
\medskip

\noindent\textbf{Construction of a frame of a simplicial complex.} Let $\Delta$ be a simplicial complex on the vertex set $\{1,\ldots,m\}$. For each integer $i$, 
let $\Gamma_{i}(\Delta)$ be the set of all $i$ dimensional faces of $\Delta$ 
and let $K^{\Gamma_{i}(\Delta)}$ denote the $K$-vector space generated by the 
basis elements $e_{\sigma}$, for $\sigma\in \Gamma_{i}(\Delta)$. The 
\textit{chain complex} of $\Delta$ over $K$ is the complex 
$\mathcal{C}(\Delta)_{\centerdot}$,
$$0\longrightarrow K^{\Gamma_{m-1}(\Delta)}\stackrel{\partial_{m-1}}{\longrightarrow}\cdots K^{\Gamma_{i}(\Delta)}\stackrel{\partial_{i}}{\longrightarrow} K^{\Gamma_{i-1}(\Delta)}\stackrel{\partial_{i-1}}{\longrightarrow} \cdots \stackrel{\partial_{0}}{\longrightarrow} K^{\Gamma_{-1}(\Delta)}\longrightarrow 0.$$
The boundary maps $\partial_{i}$ are defined as 
$$\partial_{i}(e_{\sigma})=\sum_{j\in\sigma}sign(j,\sigma)e_{\sigma-j},\quad 0\leq i\leq m-1;$$
such that $sign(j,\sigma)=(-1)^{t-1}$, where $j$ is the $t$-th element of the set 
$\sigma\subset \{1,\ldots,m\}$ written in increasing order. If $i<-1$ or 
$i>m-1$, then $K^{\Gamma_{i}(\Delta)}=0$ and  $\partial_{i}=0$. We see that 
the chain complex of a simplicial complex is an $m$-frame.
\medskip

If $\Delta$ is a standard simplicial complex of dimension $m-1$ 
on the vertex set $\{1, \ldots, m\}$, then 
$\mathcal{C}(\Delta)_{\centerdot}$ is acyclic and 
$K^{\Gamma_{i}(\Delta)}=K^{\binom{m}{i}}$. Therefore, 
the $m$-frame of a standard simplicial complex of 
dimension $m-1$ is nothing but the Koszul complex.
\medskip

\begin{definition}
Let $M=\{m_{1},m_{2},\cdots,m_{r}\}$ be a set of monomials in $R$. 
An $M-$complex $\mathbb{G}_{\centerdot}$ is a multigraded complex of 
finitely generated free multigraded $R$ modules with 
differential $\Delta$ and a fixed multihomogeneous 
basis with multidegrees that satisfy 
\begin{enumerate}[(i)]
\item $\mathbb{G}_{i}=0$ \, for \, $i<0$;
\item $\mathbb{G}_{0}=R$;
\item $\mathbb{G}_{1}=R(m_{1})\oplus \cdots\oplus R(m_{r})$;
\item $\Delta(w_{j})=m_{j}$ for each basis vector $w_{j}\in\mathbb{G}_{1}$.
\end{enumerate}
\end{definition}

\noindent\textbf{Construction.}  Let $\mathbb{U}_{\centerdot}$ be an 
$r-$frame and $M=\{m_{1},m_{2},\cdots,m_{r}\}$ be a set of monomials 
in $R$. The $M-$\textit{homogenization} of $\mathbb{U}_{\centerdot}$ 
is defined to be the complex $\mathbb{G}_{\centerdot}$ with 
$\mathbb{G}_{0}= R$ and $\mathbb{G}_{1}= R(m_{1})\oplus \cdots\oplus R(m_{r})$. 
Let $\bar{v_{1}},\cdots, \bar{v_{p}}$ and $\bar{u_{1}},\cdots, \bar{u_{q}}$ 
be the given bases of $\mathbb{U}_{i}$ and $\mathbb{U}_{i-1}$ 
respectively. Let $u_{1},u_{2},\cdots, u_{q}$ be the basis of 
$\mathbb{G}_{i-1}$ chosen in the previous step of induction. 
We take $v_{1},\cdots,v_{p}$ to be a basis of $\mathbb{G}_{i}$. 
If $$\partial(\bar{v_{j}})=\sum_{1\leq s\leq q}\alpha_{sj}\bar{u_{s}}$$
with coefficients $\alpha_{sj}\in k$, then set
$\rm{mdeg}(v_{j})= \lcm \left(mdeg(u_{s})|\alpha_{sj}\neq 0\right)$, 
where $\lcm(\emptyset)=1$. We define 
$$\mathbb{G}_{i}=\oplus_{1\leq i\leq p}R\left(\rm{mdeg}(v_{j})\right), \quad 
\Delta(v_{j})= \Sigma_{1\leq s\leq q} \alpha_{sj} \dfrac{\rm{mdeg} (v_{j})}{\rm{mdeg}(u_{s})}u_{s}.$$
\medskip

\begin{definition} Let $I\subset R$ be a monomial ideal with $G(I)=\{m_{1},\ldots,m_{r}\}$. We say that a minimal free resolution of $I$ is supported by a simplicial complex 
$\Delta$ on the vertex set $\{1,\ldots,r\}$ if it is isomorphic to the $G(I)$ homogenization of the chain complex $\mathcal{C}(\Delta)_{\centerdot}$ of  $\Delta$. 
\end{definition}
\medskip

\begin{lemma}
Let $I$ and $J$ be two monomial ideals such that $I\cap J=IJ$, $|G(I)|=r$ and $|G(J)|=s$, and whose minimal free resolutions are supported by standard simplicial complexes $\Delta_{I}$ and $\Delta_{J}$ on the vertex sets $\{1,\ldots,r\}$ and $\{r+1,\ldots,r+s\}$ respectively. Then minimal free resolution of $I+J$ is supported by the standard simplicial complex $\Delta_{I}*\Delta_{J}$. 
\end{lemma} 

\proof We have seen that if $\Delta$ is a standard simplicial complex of 
dimension $m-1$ on the vertex set $\{1, \ldots, m\}$, then the 
$m$-frame of $\Delta$ is nothing but the Koszul complex. Therefore, it follows that 
if a minimal free resolution $\mathbb{M}(I)_{\centerdot}$ of $I$ 
is supported by the standard simplicial complex $\Delta_{I}$, then the it is actually isomorphic to 
the Taylor complex $\mathbb{T}(I)_{\centerdot}$. Since  $I\cap J=IJ$, by the Theorem \ref{disjoint} we have $G(I)\cap G(J)=\emptyset$, then  $|G(I)\cup G(J)|=r+s$. Again $\Delta_{I}*\Delta_{J}$ is also standard simplicial complex of dimension $r+s-1$. Therefore $G(I)\cup G(J)$-homogenization of $\mathcal{C}(\Delta_{I}*\Delta_{J})_{\centerdot}$ is isomorphic to the Taylor complex $\mathbb{T}(I+J)_{\centerdot}$. Thus by Theorem \ref{tcomplex}, we have $$\mathbb{M}(I)_{\centerdot}\otimes\mathbb{M}(J)_{\centerdot}\cong \mathbb{T}(I)_{\centerdot}\otimes\mathbb{T}(J)_{\centerdot}\cong \mathbb{T}(I+J)_{\centerdot}$$
Again by the Corollary \ref{tenminimal}, $\mathbb{M}(I)_{\centerdot}\otimes\mathbb{M}(J)_{\centerdot}$ is minimal free resolution of $I+J$. 
Therefore minimal free resolution of $I+J$ is supported by the standard simplicial complex $\Delta_{I}*\Delta_{J}$\qed
\medskip

\begin{lemma}\label{joinframe}
Let $\Delta_{I}$ be a simplicial complex on vertex $\{1,\ldots,r\}$ and $\Delta_{J}$ be another simplicial complex on vertex $\{r+1,\ldots,r+s\}$,i.e., the vertex sets of 
$\Delta_{I}$ and $\Delta_{J}$ are disjoint. Then 
$$\mathcal{C}(\Delta_{I})_{\centerdot}\otimes \mathcal{C}(\Delta_{J})_{\centerdot}\cong \mathcal{C}(\Delta_{I}*\Delta_{J})_{\centerdot}$$ 
\end{lemma}

\proof Let $\Gamma(\Delta_{I})=\{\gamma_{1},\ldots,\gamma_{l_{1}}\}$ and 
$\Gamma(\Delta_{J})=\{\sigma_{1},\ldots,\sigma_{l_{2}}\}$. We have 
$\Gamma(\Delta_{I}*\Delta_{J})= \{\gamma_{i}\cup \sigma_{j}\mid 1\leq i\leq l_{1}, 1\leq j\leq l_{2}\}$. Now assuming $\binom{s}{t}=0$ for $s<t$, we have 
$\mid\Gamma_{p}(\Delta_{I})\mid = \sum_{i=1}^{l_{1}}\binom{\mid \gamma_{i}\mid}{p}$, for all $1\leq p\leq l_{1}$ and  $\mid\Gamma_{q}(\Delta_{J})\mid = \sum_{j=1}^{l_{2}}\binom{\mid \sigma_{j}\mid}{q}$ for all $1\leq q\leq l_{2}$. On the other hand, 
$\mid\Gamma_{t}(\Delta_{I}*\Delta_{J})\mid=\sum_{j=1}^{l_{2}}\sum_{i=1}^{l_{1}}\binom{\mid\gamma_{i}\cup \sigma_{j}\mid}{t}= \sum_{j=1}^{l_{2}}\sum_{i=1}^{l_{1}}\binom{\mid\gamma_{i}\mid+\mid\sigma_{j}\mid}{t}=\sum_{p+q=t}\sum_{j=1}^{l_{2}}\sum_{i=1}^{l_{1}}\binom{\mid\gamma_{i}\mid}{p}\binom{\mid\sigma_{j}\mid}{q}$. Therefore the map $$\theta_{t}:\left(\mathcal{C}(\Delta_{I})_{\centerdot}\otimes \mathcal{C}(\Delta_{J})_{\centerdot}\right)_{t}\longrightarrow \left(\mathcal{C}(\Delta_{I}*\Delta_{J})_{\centerdot}\right)_{t}$$ defined by $$\theta_{t}(e_{\gamma}\otimes e_{\sigma})=e_{\gamma\cup \sigma}$$ is an isomorphism, where $\gamma\in \Delta_{I}, \mid\gamma\mid=p $ and $\sigma\in \Delta_{J}, \mid\gamma\mid=q $ and $p+q=t$. We therefore have to show the following 
diagram commutes 
\medskip

\begin{center}
\begin{tikzcd}[column sep = 3cm]
(\mathcal{C}(\Delta_{I})_{\cdot}\otimes \mathcal{C}(\Delta_{J})_{\cdot})_{t} 
\arrow[r, "(\partial(\Delta_{I})\otimes\partial(\Delta_{J}))_{t}"] 
\arrow[d,"\theta_{t}"] 
& (\mathcal{C}(\Delta_{I})_{\cdot}\otimes \mathcal{C}(\Delta_{J})_{\cdot})_{t-1}
\arrow[d,"\theta_{t-1}"]\\
\mathcal{C}(\Delta_{I}*\Delta_{J})_{t}\arrow[r, "(\partial(\Delta_{I}*\Delta_{J}))_{t}"] & \mathcal{C}(\Delta_{I}*\Delta_{J})_{t-1}
\end{tikzcd}
\end{center}
\medskip

Let $\gamma\in \Delta_{I}, \mid\gamma\mid=p $ and $\sigma\in \Delta_{J}, \mid\gamma\mid=q$ 
such that $p+q=t$. Then,
\begin{eqnarray*}
{} & \theta_{t-1}\circ \left(\partial(\Delta_{I})\otimes\partial(\Delta_{J})_{t}\right)(e_{\gamma}\otimes e_{\sigma})\\
= & \theta_{t-1}((\partial(\Delta_{I})_{p}e_{\gamma})\otimes e_{\sigma}+(-1)^{p}e_{\gamma}\otimes(\partial(\Delta_{J})_{p}e_{\sigma}))\\
= & \theta_{t-1}(\sum_{j\in \gamma}sign(j,\gamma)e_{\gamma-j}\otimes e_{\sigma}+(-1)^{p}e_{\gamma}\otimes(\sum_{s\in \sigma}sign(s,\sigma)e_{\sigma-s}))\\
= & \sum_{j\in \gamma}sign(j,\gamma)e_{\gamma\cup\sigma-j}+(-1)^{p}(\sum_{s\in \sigma}sign(s,\sigma)e_{\gamma\cup\sigma-s})\\
= & \sum_{j\in \gamma\cup\sigma}sign(j,\gamma\cup\sigma)e_{\gamma\cup\sigma-j}=\partial(\Delta_{I}*\Delta_{J})_{t}\circ \theta_{t}((e_{\gamma}\otimes e_{\sigma}))\qed
\end{eqnarray*}
\medskip

\begin{theorem}\label{simpsupport}
Let $I$ and $J$ be two monomial ideals such that $I\cap J=IJ$, $|G(I)|=r$ and $|G(J)|=s$, and their minimal free resolutions are supported by the simplicial complexes $\Delta_{I}$ and $\Delta_{J}$ on the vertex sets $\{1,\ldots,r\}$ and $\{r+1,\ldots,r+s\}$ respectively. Then minimal free resolution of $I+J$ is supported by the simplicial complex $\Delta_{I}*\Delta_{J}$.
\end{theorem}

\proof We have $I\cap J=IJ$, therefore $\supp(G(I)) \cap \supp(G(J))=\emptyset$, 
by Theorem \ref{disjoint}. Again, by Lemma \ref{joinframe}, 
$$\mathcal{C}(\Delta_{I})_{\centerdot}\otimes \mathcal{C}(\Delta_{J})_{\centerdot}\cong 
\mathcal{C}(\Delta_{I}*\Delta_{J})_{\centerdot}$$
Since $\mathbb{M}(I)_{\centerdot}$ is the $G(I)$-homogenization of the complex 
$\mathcal{C}(\Delta_{I})_{\centerdot}$ and $\mathbb{M}(J)_{\centerdot}$ is the $G(J)$-homogenization of the complex $\mathcal{C}(\Delta_{J})_{\centerdot}$, we can 
proceed in the same way as lemma \ref{joinframe} to prove that $\mathbb{M}(I)_{\centerdot}\otimes \mathbb{M}(J)_{\centerdot}$ is the $G(I)\cup G(J)$--homogenization of the complex $\mathcal{C}(\Delta_{I}*\Delta_{J})_{\centerdot}$. Again by the Corollary \ref{tenminimal}, $\mathbb{M}(I)_{\centerdot}\otimes \mathbb{M}(J)_{\centerdot}$ is minimal free resolution of $I+J$. 
Therefore, minimal free resolution of $I+J$ is supported by the simplicial complex 
$\Delta_{I}*\Delta_{J}$. \qed
\bigskip

\section*{Funding Sources}
This research is supported by the research project EMR/2015/000776 
sponsored by the SERB, Government of India. 
\bigskip

\bibliographystyle{amsalpha}

\end{document}